\documentclass[10pt,a4paper]{article}
\usepackage{indentfirst}
\usepackage{amsmath,amssymb,amsfonts,amsthm,graphics,amssymb,amsmath,graphicx,graphics,amsthm,amsfonts,mathrsfs,multicol,multirow, amscd,color,enumerate,caption}
\usepackage{makeidx}
\usepackage{color}
\newtheorem{theorem}{Theorem}[section]

\newtheorem{lemma}[theorem]{Lemma}

\newtheorem{remark}[theorem]{Remark}

\def\proof{{\bf Proof.}\ }
\def\endproof{\hfill \ \ \hbox{$\sqcup$}\llap{\hbox{$\sqcap$}}}

\def \L {{\cal L}}
\def \H {{\cal H}}
\def \F {{\cal F}}

\begin{document}
\title{\Large\bf The $K_{1,2}$-structure-connectivity of graphs}
\author{Xiao Zhao$^{a}$\footnote{Corresponding author. Email: zhaoxiao05@126.com}, Haojie Zheng$^{a}$, Hengzhe Li$^{a}$\\
\small $^{a}$College of Mathematics and Information Science,\\
		\small Henan Normal University, Xinxiang 453007, P.R. China\\
        \small Email:zhaoxiao05@126.com, zhj9536@126.com, lihengzhe@htu.edu.cn\\
		}
\maketitle
\begin{abstract}

In this paper, we mainly investigate $K_{1,2}$-structure-connectivity for any connected graph. Let $G$ be a connected graph with $n$ vertices, we show that $\kappa(G; K_{1,2})$ is well-defined if $diam(G)\geq 4$, or $n\equiv 1\pmod 3$, or $G\notin \{C_{5},K_{n}\}$ when $n\equiv 2\pmod 3$, or there exist three vertices $u,v,w$ such that $N_{G}(u)\cap (N_{G}(v,w)\cup\{v,w\})=\emptyset$ when $n\equiv 0\pmod 3$.  Furthermore, if $G$ has $K_{1,2}$-structure-cut,  we prove $\kappa(G)/3\leq\kappa(G; K_{1,2})\leq\kappa(G)$.


\end{abstract}

\section{Introduction}

All the graphs considered in this paper are simple and undirected. Graph theoretical terminology and notion we use, not surprisingly, is essentially that of Bondy and
Murty \cite{Bondy}. In graph theory, connectivity is an important indicator of the reliability and fault  tolerability of a network.
Recall that a path with end vertices $x$ and $y$ is a {\it $(x,y)$-path}. Two $(x,y)$-paths $P$ and $Q$ in $G$ are {\it internally disjoint}
if they have no internal vertices in common. The {\it local connectivity} between distinct vertices $x$ and $y$ is the maximum number of pairwise internally disjoint $(x,y)$-paths, denoted $p(x, y)$; the local connectivity is undefined when $x=y$. A nontrivial graph $G$ is $k$-connected if $p(u,v)\geq k$ for any two distinct vertices $u$ and $v$.  The {\it connectivity} $\kappa(G)$ of $G$ is the maximum value of $k$ for which $G$ is $k$-connected. Thus, for a nontrivial graph $G$, $\kappa(G):=\min\{p(u, v) :  u, v\in V,  u\neq v\}$. By Menger's Theorem, the connectivity of a graph equals to the minimum size of a vertex set such that $G-S$ is disconnected or trivial. Classical connectivity has been widely studied in the literature from a combinatorial point of view. For the details and related topics of the connectivity of graphs, we refer the
reader to survey several articles \cite{c1,c2,c3,c4}.

Recently, as another variant of traditional connectivity, Lin et al. \cite{Lin} introduced structure connectivity and substructure connectivity as a way to evaluate the fault tolerance of a network from the perspective of individual vertices, as well as some special structure of the network. A set $\mathcal{F}$ of vertex-disjoint connected subgraphs of $G$ is a {\it subgraph-cut} of $G$ if $G-V(\mathcal{F})$ is disconnected or trivial. Let $H$ be a connected subgraph of $G$. Then $\mathcal{F}$ is an {\it $H$-structure-cut} if $\mathcal{F}$ is a subgraph-cut, and every element in $\mathcal{F}$ is isomorphic to $H$. The {\it $H$-structure-connectivity} of $G$, denoted by $\kappa(G; H)$, to be the minimum cardinality of all $H$-structure-cuts of $G$, that is, the minimum number of vertex-disjoint $H$s whose deletion results in a disconnected graph or trivial graph. It is easy to see that not every graph has a $H$-structure-cut (we shall study the existence of a $H$-structure-cut in a graph), that is, $H$-structure-cut has very strong requirements. The structure connectivity and substructure connectivity of some well-known networks, such as circulant graphs, hypercubes, cross cube, star graphs, wheel networks, divide-and-swap, etc \cite{Chelvam,Lin,LiLin,Feng,Pan,Zhang,Zhou}. Li et al. \cite{Li} studied the matching-connectivity (i.e. $K_{1,1}$-structure connectivity) of general graphs. In this paper, we shall study the $K_{1,2}$-structure connectivity of general graphs.

The remainder of this paper is organized as follows.  In Section 2, we summarize some useful
notations.
In Section 3, we study the existence of $K_{1,2}$-structure-cut for graph $G$ with $n$ vertices, and prove that $\kappa(G; K_{1,2})$ is well-defined if $diam(G)\geq 4$, or $n\equiv 1\pmod 3$, or $G\notin \{C_{5},K_{n}\}$ when $n\equiv 2\pmod 3$, or there exist three vertices $u,v,w$ such that $N_{G}(u)\cap (N_{G}(v,w)\cup\{v,w\})=\emptyset$ when $n\equiv 0\pmod 3$. In particular, $G\notin \{C_{5},K_{n}\}$ is a  sufficient and necessary condition for the existence of $K_{1,2}$-structure-cut when $n\equiv 2\pmod 3$.
In Section 4, we study the relationship of $K_{1,2}$-structure-connectivity $\kappa(G; K_{1,2})$ and connectivity $\kappa(G)$ of graph $G$, and show that $\kappa(G)/3\leq\kappa(G; K_{1,2})\leq\kappa(G)$.

\section{Definitions and Notations}

 Let $G$ be a connected graph. The vertex set of $G$ is denoted by $V(G)$ and the edge set by $E(G)$, respectively. Complete graphs are denoted by $K_n$, paths by $P_n$, cycles by $C_n$, where $n$ is the number of vertices. The length  of path $P$ is denoted by $\ell(P)$.

Let $H$ be a connected subgraph of $G$. For any vertex $u\in V(G)$, we say $u$ is {\it covered }by $H$ if $u\in V(H)$, and $u$ is {\it adjacent} to $H$ if $N_{G}(u)\cap V(H)\neq \emptyset$. Of course, $u$ is clearly covered by $H$ if $u\in V(H)$. Let $\F=\{F_1,F_2,\ldots,F_t\}$ be a set of connected subgraphs of $G$ which are isomorphic to $H$. For $S\subseteq V(G)$, we call $\F$
an {\it $H$-covering} of $S$ if $\F$ can cover all vertices of $S$, i.e., $S\subseteq V(\F)$.
For any vertex $v$ of $V(\F)$, we define that $d_{\F}(v)=| \{i|v\in V(F_i)\}|$. In other words, $d_{\F}(v)$ calculates the number of times $v$ occurred in $\F$.
In particular, $\F$ is {\it vertex-disjoint} if $d_{\F}(v)=1$ for any $v\in V(\F)$,  that is, $V(F_i)\cap V(F_j)=\emptyset$ for any $1\leq i,j\leq t$.

\section{The existence of $K_{1,2}$-structure-cut}

In this section, we discuss the existence of $K_{1,2}$-structure-cut on connected graphs. Firstly, we give a  sufficient condition of the existence of $K_{1,2}$-structure-cut for any connected graph. Let $H$ be a subgraph of graph $G$.

\begin{theorem}\label{Well-defined-1}
$\kappa(G; K_{1,2})$ is well-defined if $diam(G)\geq 4$. Moreover, our condition is optimal.
\end{theorem}
\proof
Let $G$ be a connected graph with $diam(G)\geq 4$ and $u, v$ be two distinct vertices of $G$ such that $dist_G(u,v)=diam(G)$. Pick a minimum vertex cut $S$ separating $u$ and $v$. Assume that $s=|S|$ and $S=\{w_1,w_2,\ldots,w_s\}$.
By Menger's Theorem, $G$ contains
$s$ internally disjoint $(u,v)$-paths, say $P_1,P_2,\ldots,P_s$.
For each $1\leq i\leq s$, pick a subpath $P'_i$ of $P_i\setminus\{u,v\}$ such that $\ell(P'_i)=2$ and $w_i\in  P'_i$, where $\ell(P'_i)$ denotes the length of $P'_i$.  Then $\{P'_1,P'_2,\ldots,P'_s\}$ is a $K_{1,2}$-structure-cut of $G$. Thus, $\kappa(G; K_{1,2})$ is well-defined.

It is easy to check that the diameter of $C_6$ is $3$, but it has no $K_{1,2}$-structure-cut. Thus our condition is optimal.
\endproof

Now, we discuss the existence of $K_{1,2}$-structure-cut on a connected graph $G$ with $n$ vertices by classifying $n$ to three cases: $n\equiv 1\pmod 3$, $n\equiv 2\pmod 3$,
and $n\equiv 0\pmod 3$.

\begin{lemma}\label{1}
Let $G$ be a connected graph with $n$ vertices and $m$ edges and $n\geq 3$. Then $\kappa(G; K_{1,2})$ is well-defined if $n\equiv 1\pmod 3$.
\end{lemma}

\proof Pick a subgraph $F_1$ isomorphic to
$K_{1,2}$ in $G$. If $G-F_1$ is disconnected, then $\{F_1\}$ is a $K_{1,2}$-structure-cut of $G$. Otherwise, pick a subgraph $F_2$ isomorphic to $K_{1,2}$ in $G-F_1$. If $G-(F_1\cup F_2)$ is disconnected, then $\{F_1,F_2\}$ is a $K_{1,2}$-structure-cut of $G$. Otherwise, repeat the above steps until $G-(F_1\cup F_2\cdots\cup F_i)$ is disconnected or trivial. Thus $\kappa(G; K_{1,2})$ is well-defined. \endproof

\begin{lemma}\label{2}
Let $G$ be a connected graph with $n$ vertices and $m$ edges and $n\geq 3$. If $n\equiv 2\pmod 3$, $\kappa(G; K_{1,2})$ is well-defined if and only if $G\notin \{C_{5}, K_{n}\}$.
\end{lemma}

\proof It is known that if $G\in \{C_5, K_n\}$, then $G$ does not have a $K_{1,2}$-structure-cut. It suffices
to prove the necessity. For sufficiency, it is easy to check $\kappa(G; K_{1,2})=1$ when $G$ is 1-connected. Thus, we can assume that $G$ is a 2-connected graph without $K_{1,2}$-structure-cut. Then $G$ is Hamilton graph when $n=5$, and let $C_5:=u_1u_2u_3u_4u_5$ be the Hamilton cycle of $G$. It is easy to see that $G\cong C_{5}$ when $n=m=5$. For $n=5$ and $m\geq 6$,  without loss of generality, we assume that $u_1u_3\in E(G)$. Then we have $\{u_2u_5,u_2u_4\}\subseteq E(G)$ since $V(G-\{u_1,u_3,u_4\})=\{u_2,u_5\}$ and $V(G-\{u_3,u_1,u_5\})=\{u_2,u_4\}$ and $G$ has no $K_{1,2}$-structure-cut. Similarly, we can get $\{u_1u_4,u_3u_5\}\subseteq E(G)$ by $\{u_2u_5,u_2u_4\}\subseteq E(G)$. Hence, $G\cong K_{5}$. Next, we will make induction on $n$ and $n\geq 8$.

For any $K_{1,2}$ of $G$, let us denote the $K_{1,2}$ by $v_1v_2v_3$ for convenience, we have $G-\{v_1,v_2,v_3\}\cong C_{5}$ or $G-\{v_1,v_2,v_3\}\cong K_{n-3}$ since $G$ has no $K_{1,2}$-structure-cut.

\emph{Case 1: $G-\{v_1,v_2,v_3\}\cong K_{n-3}$.}

Since $G$ is 2-connected, $|N_G(V(K_{n-3}))\cap \{v_1,v_2,v_3\}|\geq 2$. We conclude that $N_G(v_i)\cap V(K_{n-3})\neq \emptyset$ and $|N_G(\{v_1,v_2,v_3\})\cap V(K_{n-3})|\geq3$. Otherwise, without loss of generality, we assume that $N_G(v_1)\cap V(K_{n-3})=\emptyset$ or $|N_G(\{v_1,v_2\})\cap V(K_{n-3})|=1$. $|N_G(V(K_{n-3}))\cap \{v_1,v_2,v_3\}|\geq 2$ implies that $G-v_1$ has two vertex-disjoint $K'_{1,2}$ and $K''_{1,2}$ covering $v_2$ and $v_3$ respectively. Then $G-K'_{1,2}-K''_{1,2}$ has two components: $v_1$ and $K_{n-7}$, contradicting that  $G$ has no $K_{1,2}$-structure-cut.

If $d(v_i)<n-1$ for $1\leq i\leq 3$, then there exists a vertex $u$ such that $v_iu\notin E(G)$. If $u\in \{v_1,v_2,v_3\}$, then $V(K_{n-3})\cup (\{v_1,v_2,v_3\}\setminus\{v_i,u\})$ can be covered by $\frac{n-2}{3}$ vertex-disjoint $K_{1,2}$s by $N_G(v_j)\cap V(K_{n-3})\neq \emptyset$ for $j=1,2,3$, contradicting that $G$ has no $K_{1,2}$-structure-cut. If $u\in V(K_{n-3})$, then $V(K_{n-3}-u)\cup (\{v_1,v_2,v_3\}\setminus \{v_i\})$ can be covered by $\frac{n-2}{3}$ vertex-disjoint $K_{1,2}$s by $N_G(v_j)\cap V(K_{n-3})\neq \emptyset$ for $j=1,2,3$ and $|N_G(\{v_1,v_2,v_3\})\cap V(K_{n-3})|\geq3$, contradicting that $G$ has no $K_{1,2}$-structure-cut. Hence $G\cong K_{n}$.

\emph{Case 2: $G-\{v_1,v_2,v_3\}\cong C_{5}$.}

Since $G$ is 2-connected, $|N_G(V(C_{5}))\cap \{v_1,v_2,v_3\}|\geq 2$. We conclude that $N_G(v_i)\cap V(C_{5})\neq \emptyset$ and $|N_G(\{v_1,v_2,v_3\})\cap V(C_{5})|\geq3$. Otherwise, without loss of generality, we assume that $N_G(v_1)\cap V(C_{5})=\emptyset$ or $|N_G(\{v_1,v_2\})\cap V(C_{5})|=1$. Since $|N_G(V(C_{5}))\cap \{v_1,v_2,v_3\}|\geq 2$, we have $|N_G(\{v_2,v_3\})\cap V(C_{5})|\geq 2$ that implies $G-v_1$ has two vertex-disjoint $K'_{1,2}$ and $K''_{1,2}$ covering  $v_2$ and $v_3$ respectively and covering the only neighbor of $v_1$. Then $G-K'_{1,2}-K''_{1,2}$ has two components: $v_1$ and $V(C_5)\setminus (V(K'_{1,2})\cup V(K''_{1,2}))$, contradicting that  $G$ has no $K_{1,2}$-structure-cut.

Since $N_G(v_i)\cap V(C_{5})\neq \emptyset$ and $|N_G(\{v_1,v_2,v_3\})\cap V(C_{5})|\geq3$, there exist three different vertices $u_1,u_2,u_3$ in $N_G(\{v_1,v_2,v_3\})$ such that $\{u_1v_1,u_2v_2, u_3v_3\}\subseteq E(G)$. We have $v_1v_3\in E(G)$ since $V(C_5)\cup\{v_2\}$ can be covered by two vertex-disjoint $K_{1,2}$s. $\{u_1,u_2,u_3\}\subseteq V(C_5)$ implies $G[\{u_1,u_2,u_3\}]\neq C_3$. Without loss of generality, we assume that $u_1u_3\notin E(G)$. There is a vertex $u_4$ with $d_{C_5-u_1}(u_4)=2$ and $u_4u_3\in E(G)$. Then $K'_{1,2}:=u_1v_1v_2$ and $K''_{1,2}:=u_4u_3v_3$ are vertex-disjoint. And $G-K'_{1,2}-K''_{1,2}$ has two components: $u_2$ and $V(C_5)\setminus\{u_1,u_2,u_3,u_4\}$, contradicting that  $G$ has no $K_{1,2}$-structure-cut. \endproof

\begin{lemma}\label{3}
Let $G$ be a connected graph with $n$ vertices and $m$ edges and $n\geq 3$. If $n\equiv 0\pmod 3$ and there exist three vertices $u,v,w$ such that $N_{G}(u)\cap (N_{G}(v,w)\cup\{v,w\})=\emptyset$, $\kappa(G; K_{1,2})$ is well-defined.
\end{lemma}

\proof For $n\equiv 0\pmod 3$, let $S=\{u,v,w\}$. $N_{G}(u)\cap (N_{G}(v,w)\cup\{v,w\})=\emptyset$ implies that the number of components of $G[S])$ is $2$ or $3$. If $G-S$ has no $K_{1,2}$-structure-cut, then $V(G-S)$ can be covered by $\frac{n}{3}-1$ vertex-disjoint $K_{1,2}$s, implying that $G$ has $K_{1,2}$-structure-cut. Next, we assume that $G-S$ has $K_{1,2}$-structure-cut and let $\F$ be a maximum $K_{1,2}$-structure-cut of $G-S$. Then each  component of $G-S-V(\F)$ is  either isolated vertex or $K_2$ and $|V(G-S-V(\F))|\geq 3$. We denote all isolated vertices of $G-S-V(\F)$ by $z_1,\ldots,z_s$ and all $K_2$s of $G-S-V(\F)$ by $O_1:=x_1y_1, O_2:=x_2y_2,\ldots,O_t:=x_ty_t$.

Let $G':=G-V(\F)$. If $G'$ is not connected, then $\F$ is also a $K_{1,2}$-structure-cut of $G$. If $G'$ is 1-connected, $G'$ has a $K_{1,2}$-structure-cut $\H$ by $|V(G')|\geq 6$. Then $\F\cup \H$ is a $K_{1,2}$-structure-cut of $G$. Next we consider that $G'$ is 2-connected, and three cases should be analysed.

\begin{figure}[h]
\centering
\scalebox{1}[1]{\includegraphics{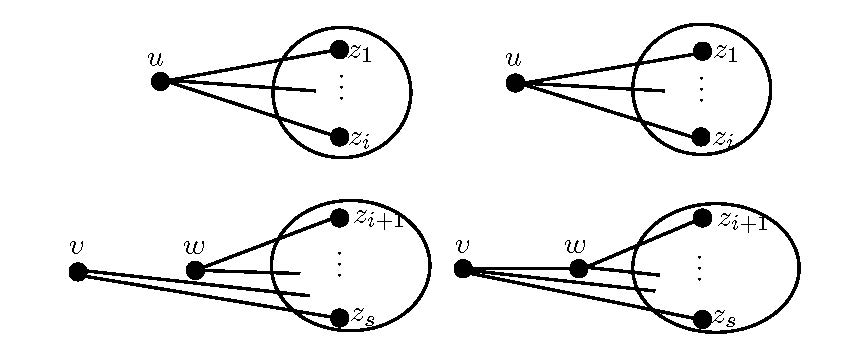}}
\caption{ The illustration of Case 1}
\label{f1}
\end{figure}

\emph{Case 1: Every component of $G-S-V(\F)$ is an  isolated vertex.} Since $N_{G}(u)\cap (N_{G}(v,w)\cup\{v,w\})=\emptyset$, we have $G'[\{u\}\cup N_{G'}(u)]$ and $G'[\{v,w\}\cup N_{G'}(v,w)]$ are two components of $G'$, as shown in Figure \ref{f1}. This contradicts that $G'$ is connected.

\emph{Case 2: Every component of $G-S-V(\F)$ is a $K_2$.}

\begin{figure}[h]
\centering
\scalebox{1}[1]{\includegraphics{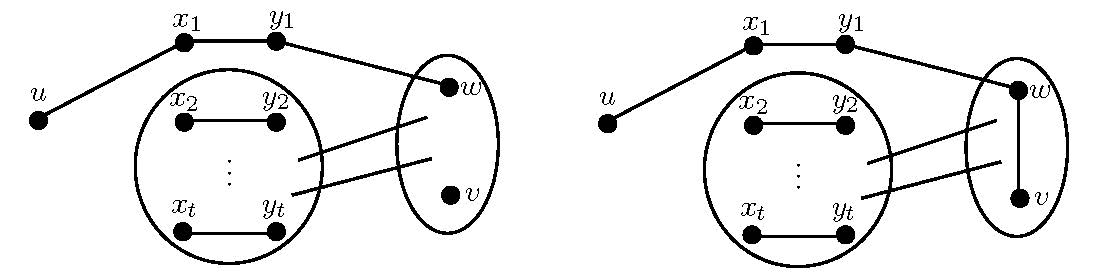}}
\caption{$\{wy_1x_1\}$ is a $K_{1,2}$-structure-cut of $G-\F$}
\label{f2}
\end{figure}

Then $n\equiv 0\pmod 3$  implies $t\geq 3$. Since $G'$ is connected, there exists $O_i$ connecting $u$ to $v$ or $w$, where $1\leq i\leq t$. Without loss of generality, we assume that $\{ux_1, wy_1\}\subseteq E(G')$. If there is no component of $G-S-V(\F)$ which is  adjacent to $u$ except for $x_1y_1$, then $G-V(\F)-\{wy_1x_1\}$ is not connected by $N_{G}(u)\cap (N_{G}(v,w)\cup\{v,w\})=\emptyset$, as shown in Figure \ref{f2}. Hence $\F\cup \{wy_1x_1\}$ is a $K_{1,2}$-structure-cut of $G$.

\begin{figure}[h]
\centering
\scalebox{1}[1]{\includegraphics{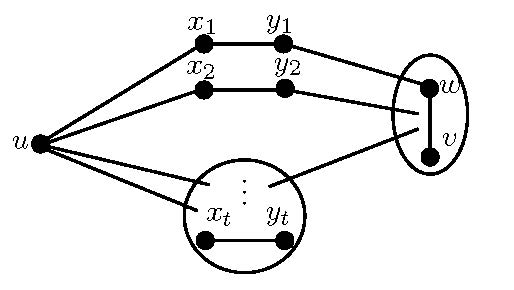}}
\caption{$\{x_1ux_2,y_1wv\}$ is  a $K_{1,2}$-structure-cut of $G-\F$}
\label{f3}
\end{figure}

If there is $O_j$ such that $2\leq j\leq t$ and $\{ux_i,uy_i\}\cap E(G)\neq \emptyset$, without loss of generality, we assume that $ux_2\in E(G)$. Then we have $G-V(\F)-\{x_1,u,x_2,y_1,w,v\}$ is not connected and $\F\cup\{x_1ux_2,y_1wv\}$ is  a $K_{1,2}$-structure-cut of $G$ when $vw\in E(G)$, as shown in Figure \ref{f3}. Next, we assume that $\{ux_1,ux_2\}\subseteq E(G')$ and $vw\notin E(G)$.

\begin{figure}[h]
\centering
\scalebox{1}[1]{\includegraphics{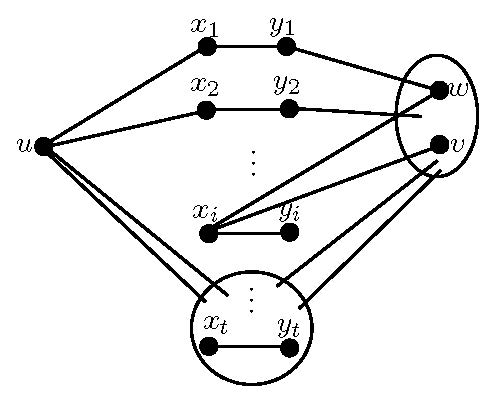}}
\caption{$\{vx_iw,x_1ux_2\}$ is  a $K_{1,2}$-structure-cut of $G-\F$}
\label{f4}
\end{figure}

If there is $x_i$ ($y_j$) such that $\{vx_i,wx_i\}\subseteq E(G')$ (resp. $\{vy_j,wy_j\}\subseteq E(G')$), where $3\leq i\leq t$ (resp. $3\leq j\leq t$). Then $G-V(\F)-\{v,x_i,w,x_1,u,x_2\}$ (resp. $G-V(\F)-\{v,y_j,w,x_1,u,x_2\}$) is not connected since $m\geq 3$. Hence $\F\cup\{vx_iw,x_1ux_2\}$ (resp. $\F\cup\{vy_jw,x_1ux_2\}$) is  a $K_{1,2}$-structure-cut of $G$, as shown in Figure \ref{f4}.

\begin{figure}[h]
\centering
\scalebox{1}[1]{\includegraphics{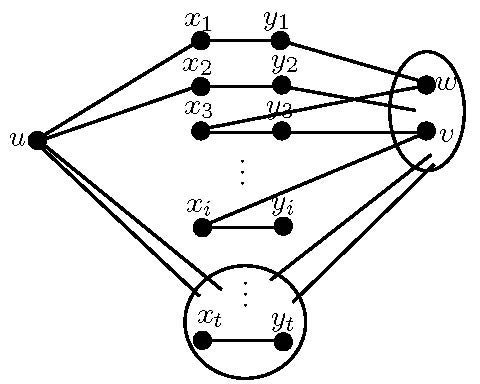}}
\caption{$\{x_1y_1w,x_ivy_3\}$ is a $K_{1,2}$-structure-cut of $G-\F$}
\label{f5}
\end{figure}

If any vertex of $\{x_3,\ldots,x_t,y_1,\ldots,y_t\}$ is only adjacent to one vertex of $\{u,v,w\}$, then we conclude that there exits $O_i:=x_iy_i$ ($3\leq i\leq t$) such that $\{wx_i,vy_i\}\subseteq E(G')$ or $\{wy_i,vx_i\}\subseteq E(G')$ (without loss of generality, we assume that $\{wx_3,vy_3\}\subseteq E(G')$). Otherwise, $G'-u$ is not connected since $vw\notin E(G)$,  contradicting that $G'$ is 2-connected. Furthermore, there exist $x_i$ or $y_j$ such that $x_iv\in E(G')$ or $y_jv
\in E(G')$ by $d_{G'}(v)\geq 2$, where $4\leq i\leq t$ and $2\leq j\leq t$ and $j\neq 3$, as shown in Figure \ref{f5}. Then $G-V(\F)-\{x_1,y_1,w,y_j,v,y_3\}$ or $G-V(\F)-\{x_1,y_1,w, x_i,v,y_3\}$ is not connected. Hence $\F\cup \{x_1y_1w,y_jvy_3\}$ or $\F\cup \{x_1y_1w, x_ivy_3\}$ is a $K_{1,2}$-structure-cut of $G$.

\emph{Case 3: There exist not only isolated vertex but also $K_2$ in $G-S-\F$.}

We conclude that $u$ should be connected to $v$ or $w$ by a $K_2$ of $G-S-\F$  since $G'$ is connected and $N_{G}(u)\cap (N_{G}(v,w)\cup\{v,w\})=\emptyset$. Without loss of generality, we assume that $ux_1,wy_1\in E(G')$.

\begin{figure}[h]
\centering
\scalebox{1}[1]{\includegraphics{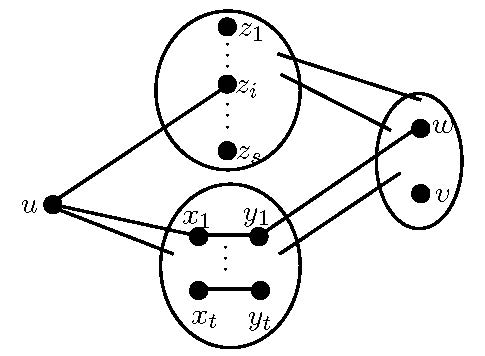}}
\caption{$\{ux_1y_1\}$ is a $K_{1,2}$-structure-cut of $G-\F$}
\label{f6}
\end{figure}

If there is an isolated vertex $z_i$ ($1\leq i\leq s$) such that $uz_i\in E(G')$, then $z_iw,z_iv\notin E(G')$ and $z_i$ is also an isolated vertex of $G-V(\F)-\{u,x_1,y_1\}$ and  $|V(G-V(\F)-\{u,x_1,y_1\})|\geq 2$, as shown in Figure \ref{f6}. Hence $\F\cup\{ux_1y_1\}$ is a $K_{1,2}$-structure-cut of $G$.

\begin{figure}[h]
\centering
\scalebox{1}[1]{\includegraphics{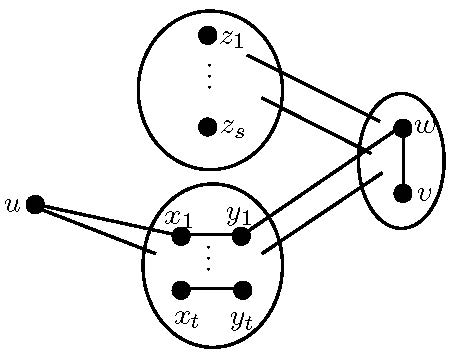}}
\caption{$\{vwy_1\}$ is a $K_{1,2}$-structure-cut of $\F$}
\label{f7}
\end{figure}

If all isolated vertices are not adjacent to $u$ and $vw\in E(G')$, as shown in Figure \ref{f7}, then those isolated vertices also are isolated vertices of $G-V(\F)-\{v,w,y_1\}$. Hence $\F\cup\{vwy_1\}$ is a $K_{1,2}$-structure-cut of $G$. Next, we assume that all isolated vertices are adjacent to $v$ or $w$ and $vw\notin E(G')$.

\begin{figure}[h]
\centering
\scalebox{1}[1]{\includegraphics{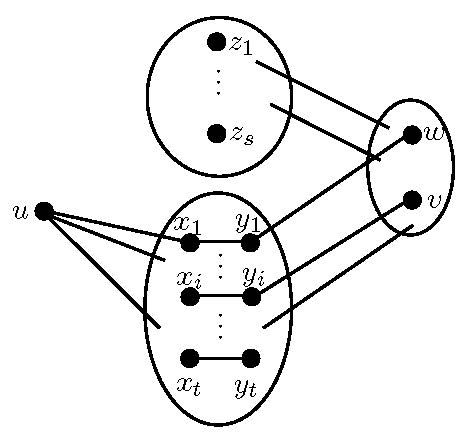}}
\caption{$\{wy_1x_1,vy_ix_i\}$ is a $K_{1,2}$-structure-cut of $G-\F$}
\label{f8}
\end{figure}

If $vy_1\in E(G')$, then $G-V(\F)-\{v,y_1,w\}$ is not connected and $\F\cup\{vy_1w\}$ is a $K_{1,2}$-structure-cut of $G$.  If there is a $O_i$ $(2\leq i\leq t)$ such that $vx_i\in E(G)$ ($vy_i\in E(G)$), as shown in Figure \ref{f8}, then $G-V(\F)-\{w,y_1,x_1,v,x_i,y_i\}$ (resp. $G-V(\F)-\{w,y_1,x_1,v,y_i,x_i\}$) is not connected. Hence $\F\cup\{wy_1x_1,vx_iy_i\}$ (resp. $F\cup\{wy_1x_1,vy_ix_i\}$) is a $K_{1,2}$-structure-cut of $G$.

\begin{figure}[h]
\centering
\scalebox{1}[1]{\includegraphics{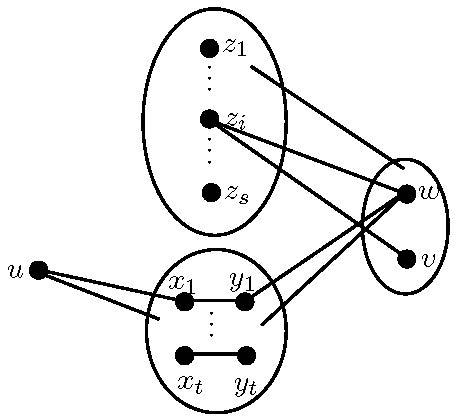}}
\caption{$\{wy_1x_1\}$ is a $K_{1,2}$-structure-cut of $G-\F$}
\label{f9}
\end{figure}

If neither of both situations  happens, there is a isolated vertex $z_i$ adjoined $v$ and $w$ since $G'$ is connected, as shown in Figure \ref{f9}. Then $G-V(\F)-\{w,y_1,x_1\}$ is not connected and $\F\cup\{wy_1x_1\}$ is a $K_{1,2}$-structure-cut of $G$.\endproof

\begin{figure}[h]
\centering
\scalebox{1}[1]{\includegraphics{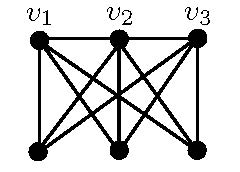}}
\caption{An example of Remark \ref{r}}
\label{e1}
\end{figure}

\begin{remark}\label{r}
The converse of the above theorem is not necessarily true. The condition that  $N_{G}(u)\cap (N_{G}(v,w)\cup\{v,w\})=\emptyset$  implies $dist_{G}(u,v)\geq 3$ and $d_{G}(u,v)\geq 3$. The converse of the above theorem is not necessarily true. Let $G$ be the graph shown in Figure \ref{e1}, $G$ has a $K_{1,2}$-structure-cut $\{v_1v_2v_3\}$. But the distance of  any two distinct vertices of $G$  is at most two.
\end{remark}

By Lemmas \ref{1}, \ref{2} and \ref{3},  we get the following result.

\begin{theorem}\label{Well-defined-2}
Let $G$ be a connected graph with $n$ vertices and $m$ edges and $n\geq 3$.

(1) If $n\equiv 1\pmod 3$, $\kappa(G; K_{1,2})$ is well-defined.

(2) If $n\equiv 2\pmod 3$, $\kappa(G; K_{1,2})$ is well-defined if and only if $G\notin \{C_{5}, K_{n}\}$.

(3) If $n\equiv 0\pmod 3$ and there exist three vertices $u,v,w$ such that $N_{G}(u)\cap (N_{G}(v,w)\cup\{v,w\})=\emptyset$, $\kappa(G; K_{1,2})$ is well-defined.
\end{theorem}

\section{$K_{1,2}$-structure-connectivity}

In this section, our main objective is to study the relationship between $\kappa(G)$ and $\kappa(G; K_{1,2})$ and show our main result: $\kappa(G)/3\leq\kappa(G; K_{1,2})\leq \kappa(G)$. As a first step towards
the proof of this statement, we show that any minimum vertex cut of $G$ has a vertex-disjoint $K_{1,2}$-covering. We exploit the property
in the proof of the assertion stated above.

\begin{lemma}\label{covering}
Let $G$ be a connected graph with $n$ vertices and $n\geq 4$ and $X$ be a vertex cut of $G$ with $|X|=\kappa(G)$. Then there is a vertex-disjoint $K_{1,2}$-covering $\F$ of $X$ such that $|\F|\leq |X|$ in $G$.
\end{lemma}

\proof It is trivial when $\kappa(G)=1$. We may thus assume that $\kappa(G)=k\geq2$. Let $X=\{v_1,v_2,\ldots,v_k\}$. Clearly, $G-X$ consists some connected components. We divide these connected components into two nonempty parts and denote their vertex sets  by $S_1$ and $S_2$, respectively.

Since $\kappa(G)\geq 2$, we have for any vertex $v\in V(G)$, $d_G(v)\geq k\geq2$. It implies that each vertex of $X$ can be covered by a $K_{1,2}$. Now we are going to prove that there exists  a vertex-disjoint $K_{1,2}$-covering of $X$.


\noindent {\bf Claim 1}: $\forall v_i\in X$ $(i=1,2,\ldots,k)$, $N_{G}(v_i)\cap S_1\neq\emptyset$ and $N_{G}(v_i)\cap S_2\neq\emptyset$.

This is obviously valid, otherwise it contradicts that $X$ is a minimal vertex cut set.

\begin{figure}[h]
\centering
\scalebox{0.9}[0.9]{\includegraphics{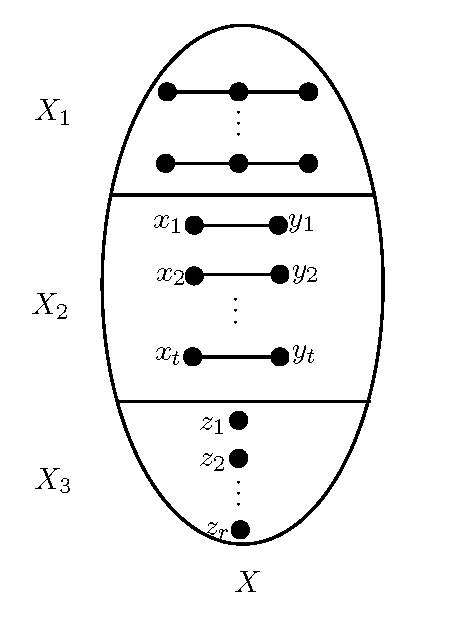}}
\caption{Three parts of $X$}
\label{t1}
\end{figure}

Now we divide $X$  into three parts by the following operations, as shown in Figure \ref{t1}:

\emph{Step 1}: Find a maximum set of vertex-disjoint $K_{1,2}$s $\F_1$ of $G[X]$, and let $X_1:=V(\F_1)$;

\emph{Step 2}: Find a maximum set of vertex-disjoint $K_{1,1}$s (i.e. matching) $M$ of $G[X\setminus X_1]$, and let $X_2:=V(M)$ and denote all elements of $M$ by $x_1y_1,x_2y_2,\ldots,x_ty_t$;

\emph{Step 3}: Let $X_3:=X\setminus(X_1\cup X_2)$ and denote all vertices of $X_3$ by $z_1,z_2,\ldots,z_r$.

%
Clearly, $X_3$ is a stable set and the vertex of $X_3$ are not adjacent to  the vertex of $X_2$. Moreover, the vertices of different elements of $M$ are not adjacent.

\noindent {\bf Claim 2}: There exists a vertex-disjoint $K_{1,2}$-covering $\F_2$ of $X_2$ such that $x_iy_i$ is a subgraph of one element of $\F_2$ and $V(\F_2)\cap X_1=\emptyset$ for $i\in \{1,\ldots,t\}$.

By Claim 1, there is a $K_{1,2}$-covering $\F_2$ of $X_2$ such that  $\forall i\in \{1,\ldots,t\}$, $x_iy_i$ is an subgraph of one element of $\F_2$ and $V(\F_2)\cap X_1=\emptyset$ and $\F_2\setminus X_2$ is as large as possible. If $|\F_2\setminus X_2|=t$, then Claim 2 holds. Otherwise, without loss of generality, we assume that $\{wx_1y_1,wx_2y_2\}\subseteq\F_2$. Then we claim that $N_G(x_1)\setminus (X_1\cup V(\F_2))\neq \emptyset$. Otherwise, we have $N_G(x_1)\subseteq X_1\cup V(\F_2)$. We denote the number of $K_{1,2}$s adjacent to $x_1$ in $\F_1$ by $h$ and the number of $K_{1,2}$s adjacent to $x_1$ in $\F_2$ except for $wx_1y_1$ by $l$. On one hand, $x_1$ is adjacent to at most three vertices of each $K_{1,2}$ in $\F_1$ and at most one vertex of each $K_{1,2}$ in $\F_2$, implying that $1+3h+l\geq d_G(x_1)$. On the other hand, $d_G(x_1)\geq k\geq 2+3h+2l.$ Hence $$1+3h+l\geq d_G(x_1)\geq k\geq 2+3h+2l.$$
Obviously this inequality does not hold.

Now we can replace $wx_1y_1$ by $w'x_1y_1$ in $\F_2$, where $w'\in N_G(x_1)\setminus (X_1\cup V(\F_2))$. Obviously, $\{w',x_1,y_1\}\cap \{w,x_2,y_2\}=\emptyset$ which contradicts that $\F_2\setminus X_2$ is maximum. Hence Claim 2 holds.

\noindent {\bf Claim 3}: There exists a $K_{1,2}$-covering $\H$ of $X_3$ such that $\forall v\in V(\H)$, $d_{\H}(v)\leq 2$ and $V(\H)\cap (X_1\cup X_2)=\emptyset$.

By claim 1, we have there exists a $K_{1,2}$-covering $\H':=\{K^i_{1,2}=u_iz_iv_i|u_i\in S_1, v_i\in S_2,z_i\in X_3,i\in\{1,\ldots,r\}\}$ of $X_3$. Note that $d_{\H'}(z_i)=1$ for $1\leq i\leq r$ and $V(\H')\cap (X_1\cup X_2)=\emptyset$.  If there is a vertex $v\in V(\H')$ with $d_{\H'}(v)\geq3$. Without loss of generality, we assume that $v\in K^1_{1,2}\cap K^2_{1,2}\cap K^3_{1,2}$, i.e., $v=v_1=v_2=v_3$.

We claim that $N_G(z_1)\setminus (X_1\cup V(\F_2)\cup V(\H'))\neq\emptyset$ or there is a vertex $u\in N_G(z_1)\cap V(\H')$, $d_{\H'}(u)=1$. Otherwise we have $N_G(z_1)\subseteq X_1\cup V(\F_2)\cup V(\H')$ and for any vertex $u\in N_G(z_1)\cap V(\H')$, $d_{\H'}(u)\geq 2$. We denote the number of $K_{1,2}$s adjacent to $z_1$ in $\F_1$ by $h$ and the number of $K_{1,2}$s adjacent to $z_1$ in $\F_2$ by $p$ and the number of $K_{1,2}$s adjacent to $z_1$ in $\H'$ except for $K^1_{1,2},K^2_{1,2}$ and $K^3_{1,2}$ by $q$. Note that $z_1$ is adjacent to at most two vertices of each $K_{1,2}$ in $\H'$, and $u$ is covered by at least two $K_{1,2}$s of $\H'$ for any vertex $u\in N_G(z_1)\cap V(\H')$, implying that $z_1$ is adjacent to at most $q$ vertices in $\H'\setminus\{K^1_{1,2},K^2_{1,2},K^3_{1,2}\}$.
Since $d_G(z_1)\geq k$, it follows that
$$2+3h+p+q \geq d_G(z_1)\geq k\geq 3+3h+2p+q$$
which is not established.

Now we can obtain a new $K_{1,2}$-covering of $X_3$ by replacing $v_1$ by $v'_1$ in $\H'$, where $v'_1\in N_G(z_1)\setminus (X_1\cup V(\F_2)\cup V(\H'))$ or $v'_1$ is a vertex of $N_G(z_1)\cap V(\H')$ with $d_{\H'}(v_1)=1$. Repeat the same procedure  until a $K_{1,2}$-covering  $\H$ of $X_3$ with $d_{\H}(u)\leq 2$ for any $u\in \H$ is obtained. Hence Claim 3 holds.

\begin{figure}[h]
\centering
\scalebox{1}[1]{\includegraphics{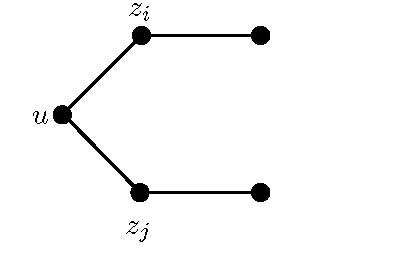}}
\caption{$z_iuz_j$ is a $K_{1,2}$-covering of $z_i$ and $z_j$}
\label{t2}
\end{figure}

Next we will obtain a vertex-disjoint $K_{1,2}$-covering $\H_1$ of $X_3$ with $V(\H_1)\cap X_1=\emptyset$ by performing the following operations in $\H$. For convenience, we still use the label of $\H'$ in $\H$, and denote the $K_{1,2}$ covering $z_i$ by $K^i_{1,2}$ for $i\in \{1,\ldots,r\}$. Let $H$ be the graph with the vertex set $W=\{v_1,v_2,\ldots,v_r\}$ and two vertices $v_i$ and $v_j$ are adjacent if and only if $V(K^i_{1,2})\cap V(K^j_{1,2})\neq\emptyset$. Let $M_1$ be a maximum matching of $H$.  If $v_iv_j\in M_1$, then we have that there exist a vertex $u\in V(K^i_{1,2})\cap V(K^j_{1,2})$, as shown in Figure \ref{t2}. Obviously $z_iuz_j$ is a new $K_{1,2}$ covering $z_i$ and $z_j$. Then we replace $K^i_{1,2}$ and $K^j_{1,2}$ by $z_iuz_j$ in $\H$. We can get a new $K_{1,2}$-covering $\H_1$ of $X_3$ by doing this for each edge of $M_1$. By Claim 3, $\H_1$ is a set of vertex-disjoint $K_{1,2}$s.

\begin{figure}[h]
\centering
\scalebox{1}[1]{\includegraphics{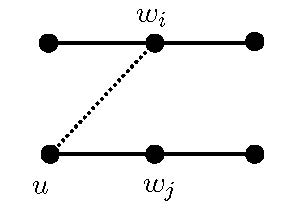}}
\caption{$w_iuw_j$ is a $K_{1,2}$-covering of $w_i$ and $w_j$}
\label{t3}
\end{figure}

Further, a better $K_{1,2}$-covering of $X_3$ can be found by using the following approach. Let $\H'_1$ be the set of $K_{1,2}$s belonging to $\H_1$ and covering exactly one vertex of $X_3$. For convenience of expression, we relabel all elements of $\H'_1$ with $H_1,H_2,\ldots,H_{|\H'_1|}$, and relabel all vertices of $V(H_i)\cap(X_3)$ with $w_1,w_2,\ldots,w_{|\H'_1|}$ temporarily. Next let $H'$ be the graph with the vertex set $W'=\{v'_1,v'_2,\ldots,v'_{|\H'_1|}\}$ and two vertices $v'_i$ and $v'_j$ are adjacent if and only if $N_G(w_i)\cap V(H_j)\neq\emptyset$ or $N_G(w_j)\cap V(H_i)\neq\emptyset$. Let $M_2$ be a maximum matching of $H'$.  If $v_iv_j\in M_2$, then we have that there exist a vertex $u\in N_G(w_i)\cap V(H_j)$ or $u\in N_G(w_j)\cap V(H_i)$.
Without loss of generality, we assume that $u\in N_G(w_i)\cap V(H_j)$, as shown in Figure \ref{t3}. Obviously $w_iuw_j$ is a new $K_{1,2}$ covering $w_i$ and $w_j$. Then we  replace $H_i$ and $H_j$ by $w_iuw_j$ in $\H'_1$. We can get a new $K_{1,2}$-covering $\H_2$ of $X_3$ by doing this for each edge of $M_2$.

According to the construction of $\H_2$, we can obtain the following result.

\noindent {\bf Claim 4}: Let $\H'_2$ be the set of $K_{1,2}$s belonging to $\H_2$ and covering exactly one vertex of $X_3$, and let $X'_3=V(\H'_2)\cap X_3$. Then we have for any vertex $z\in X'_3$, $z$ is not  adjacent to any vertex of $V(\H'_2)\setminus X'_3$ except for two vertices belonging to the $K_{1,2}$ of $\H_2$ which covers $z$.



\noindent {\bf Claim 5}: There exists a $K_{1,2}$-covering $\L_1$ of $X'_3$  with the properties of $\H'_2$ such that $K'_{1,2}\cap K''_{1,2}=\emptyset$ for any $K'_{1,2}\in \F_2$ and $K''_{1,2}\in \L_1$.

If $\H'_2$ satisfies the requirements of $\L_1$, then let $\L_1=\H'_2$  and Claim 5 holds. Otherwise, there are $K'_{1,2}\in \F_2$ and $K''_{1,2}\in \H'_2$ such that $K'_{1,2}\cap K''_{1,2}=\{w\}\neq\emptyset$. Suppose $K''_{1,2}$ covers $z$, we claim $N_G(z)\setminus V(\F_1\cup \F_2\cup \H_2)\neq \emptyset$. Otherwise, $N_G(z)\subseteq V(\F_1\cup \F_2\cup \H_2)$.  Then we denote the number of $K_{1,2}$s adjacent to $z$ of $\F_1$ by $h$,  and the number of $K_{1,2}$s adjacent to $z$ of $\F_2$ by $p$ ($p\geq1$) and the number of $K_{1,2}$s adjacent to $z$ of $\H_2\setminus\H'_2$ by $l$. By Claim 4 and $d_G(z)\geq k$, we have
$$1+3h+p+l\geq d_G(z)\geq k\geq 1+3h+2p+2l$$
which is not hold since $p\geq1$.

Now we can replace $w$ of $K''_{1,2}$ by $w'$  in $\H'_2$, where $w'\in N_G(z)\setminus V(\F_1\cup \F_2\cup \H_2)$. Hence Claim 5 holds.

Let $\L_2=\H_2\setminus\H'_2$ and $\L=\L_1\cup L_2$. Obviously, $\L$ is a vertex-disjoint $K_{1,2}$-covering of $X_3$ satisfying Claims 4 and 5, and $\L_1$ be the set of $K_{1,2}$s belonging to $\L$ and covering exactly one vertex of $X_3$, and $\L_2$ be the set of $K_{1,2}$s belonging to $\L$ and covering exactly two vertices of $X_3$.

\noindent {\bf Claim 6}: There exists a $K_{1,2}$-covering $\F_3$ of $X_3$  with the properties of $\L$ such that $K'_{1,2}\cap K''_{1,2}=\emptyset$ for any $K'_{1,2}\in \F_2$ and $K''_{1,2}\in \F_3$.

If $\L$ satisfies the requirements of $\F_3$, then let $\F_3=\L$  and Claim 6 holds. Otherwise, there are $K'_{1,2}\in \F_2$ and $K''_{1,2}\in \L$ such that $K'_{1,2}\cap K''_{1,2}\neq\emptyset$. By Claim 5, $K''_{1,2}\in \L_2$, implying that $K''_{1,2}$ covers exactly two vertices of $X_3$. Assume that $K'_{1,2}=wx_1y_1$. Then we claim $N_G(y_1)\setminus V(\F_1\cup \F_2\cup \L)\neq\emptyset$.  Otherwise $N_G(y_1)\subseteq V(\F_1\cup F_2\cup \L)$. We denote the number of $K_{1,2}$s adjacent to $y_1$ of $\F_1$ by $h$,  and the number of $K_{1,2}$s adjacent to $y_1$ of $\F_2$  by $p$  except of $K'_{1,2}$, and the number of $K_{1,2}$s adjacent to $y_1$ of $\L_2$ by $l$ except of $K''_{1,2}$, and the number of $K_{1,2}$s adjacent to $y_1$ of $\L_1$ by $q$. By $d_G(y_1)\geq k$, it follows that
$$2+3h+p+l+2q\geq d_G(y_1)\geq k\geq 4+3h+2p+2l+q.$$

If $q=0$, then it is easy to get the contradiction. Then we can replace $w$ of $K'_{1,2}$ by $w'$  in $\L$, where $w'\in N_G(y_1)\setminus V(\F_1\cup \F_2\cup \L)$.

If $q\neq 0$, there is a $K'''_{1,2}$ in $\L_1$ which covers $z'$ and  satisfies $V(K'''_{1,2})\cap N_G(y_1)\neq \emptyset$. Then we conclude that $N_G(z') \setminus V(\F_1\cup \F_2\cup \L)\neq\emptyset$. Otherwise $N_G(z')\subseteq V(\F_1\cup \F_2\cup \L)$. We denote the number of $K_{1,2}$s adjacent to $z'$ of $\F_1$ by $h'$,
 and the number of $K_{1,2}$s adjacent to $z'$ of $\F_2$ by $p'$ except of $K'_{1,2}$,
and the number of $K_{1,2}$s adjacent to $z'$ of $\L_2$ by $l'$ except of $K''_{1,2}$. By Claims 4 and 5 and $d_G(z')\geq k$, it follows that
$$3+3h'+p'+l'\geq d_G(z')\geq k\geq 5+3h'+2p'+2l',$$  a contradiction. Then we can replace the vertex of $V(K'''_{1,2})\cap N_G(y_1)$ by the vertex of $N_G(z') \setminus V(\F_1\cup \F_2\cup \L)$ in $\L$.  Repeat the same procedure  until $q=0$. Hence Claim 6 holds.

Given the above, $\F=\F_1\cup\F_2\cup\F_3$ is a vertex-disjoint $K_{1,2}$-covering of $X$ in $G$ and $|\F|\leq k$.\endproof

Now we are ready to prove the following theorem.

\begin{theorem}\label{result}
Let $G$ be a connected graph with $n$ vertices and $n\geq 4$. If $G$ has a $K_{1,2}$-structure-cut, then we have $\kappa(G)/3\leq\kappa(G; K_{1,2})\leq \kappa(G)$. Moreover, the bound is sharp.
\end{theorem}

\proof It is trivial when $\kappa(G)=1$ or $n=4$. We may thus assume that $\kappa(G)=k\geq2$. The vertex set of  a $K_{1,2}$-structure-cut of $G$  is  also a vertex cut of $G$,  implying that $\kappa(G)/3\leq\kappa(G; K_{1,2})$. Moreover, an example of graph with $\kappa(G; K_{1,2})=\kappa(G)/3$ is given in Figure \ref{e1}.

Let $X$ be a vertex cut of $G$ with $|X|=k$. By Lemma \ref{covering}, there is a vertex-disjoint $K_{1,2}$-covering $\F$ of $X$ and $|\F|\leq k$. If $G-V(\F)$ is disconnected, then $\F$ is a $K_{1,2}$-structure-cut of $G$, implying $\kappa(G; K_{1,2})\leq \kappa(G)$.

Then we  consider  $H=G-V(\F)$ as a connected graph. For simplicity, $S_1$ and $S_2$ of Lemma \ref{covering} are still used in the following discussion. Suppose that $\kappa(G; K_{1,2})=s> \kappa(G)=k$ and $V(H)\subseteq S_2$. Since $s>k$ and $|\F|\leq k$, then $|V(G)|\geq 3s+1\geq 3k+4$, thus $|V(H)|\geq 4$.

\begin{figure}[h]
\centering
\scalebox{0.9}[0.9]{\includegraphics{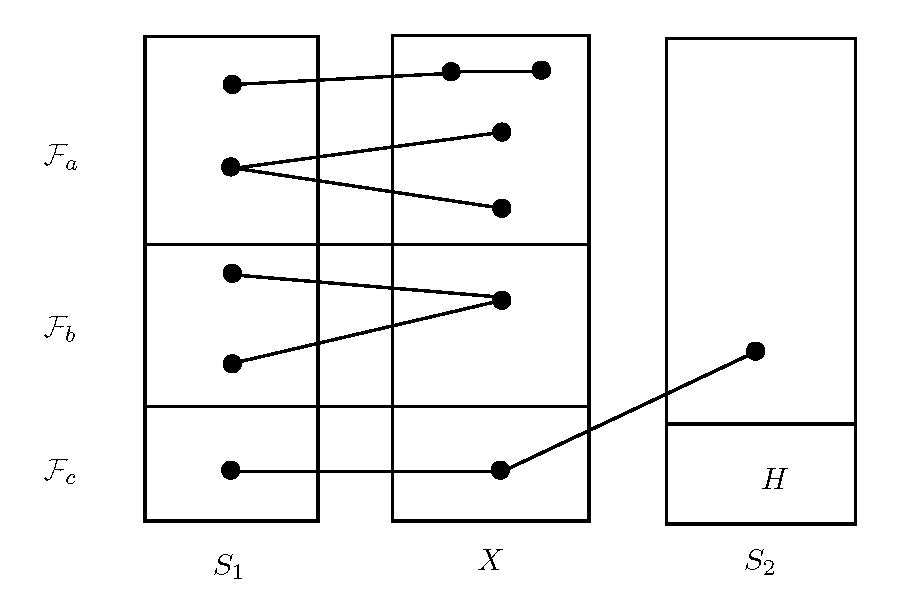}}
\caption{$\F_a$, $\F_b$ and $\F_c$}
\label{t4}
\end{figure}

Let $\F_a$ be the set of $K_{1,2}$s that covers exactly two vertices of $X$ and one vertex of $S_1$, and $\F_b$ be the set of $K_{1,2}$s that covers exactly two vertices of $S_1$, and $\F_c$ be the set of $K_{1,2}$s that covers exactly one vertex of $X$ and one vertex of $S_1$. Obviously, at least one of $\F_a$, $\F_b$ and $\F_c$ is not nonempty. As shown in Figure \ref{t4}.

\emph{Case 1: $\F_a\neq\emptyset$. }

In this case, $|\F|<k$. For any $F$ of $\F_{a}$, let $F\cap X=\{w_1,w_2\}$. If $\kappa(H)=1$, there is a $K'_{1,2}$ such that $H-\{K'_{1,2}\}$ is not connected or trival. Then we have $\F^*=\F\cup\{K'_{1,2}\}$ is a $K_{1,2}$-structure-cut of $G$ and $|\F^*|=|\F|+1\leq k$, a contradiction. Next we consider that $H$ is 2-connected. If $N_G(w_1, w_2)\cap V(H)= \emptyset$, it is easy to see that $\F^*= \F\backslash\{F\}$ is a $K_{1,2}$-structure-cut of $G$ and $|\F^*|< k$ since $G-\F^*$ is disconnected and $K_{1,2}$ is a component of $G-\F^*$. This contradicts our hypothesis. If $N_G(\{w_1, w_2\})\cap V(H)\neq \emptyset$, two cases should be analysed.

\emph{Subcase 1.1: $w_1w_2\in E(G)$.}

Pick $u'\in N_G(\{w_1, w_2)\}\cap V(H)$. Then $\F^*= (\F\backslash V(F))\cup \{u'w_1w_2\}$ is a $K_{1,2}$-structure-cut of $G$ and $|\F^*|\leq k$ since $G-\F^*$ is disconnected and $u$ is a component of $G-\F^*$. This contradicts our hypothesis.

\emph{Subcase 1.2: $w_1w_2\notin E(G)$.}

\begin{figure}[h]
\centering
\scalebox{0.9}[0.9]{\includegraphics{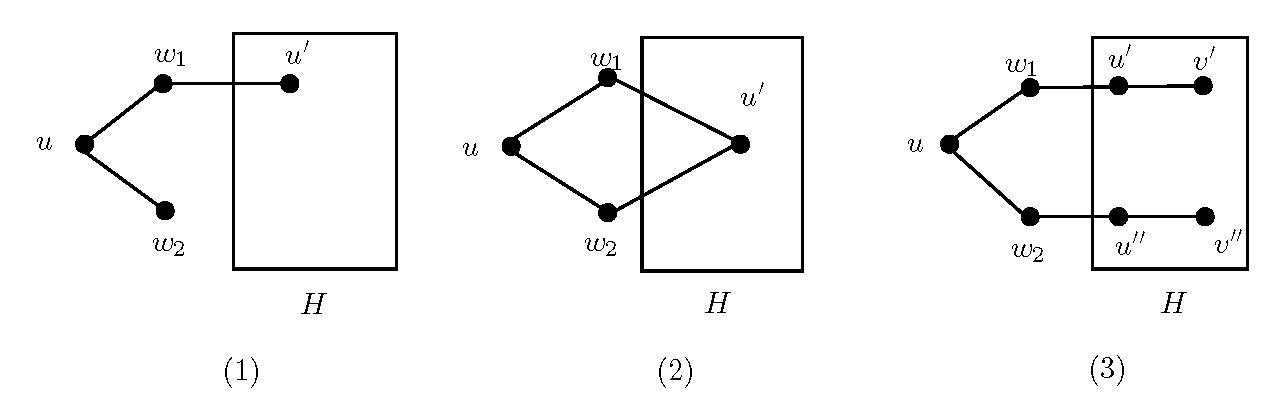}}
\caption{The illustration of Subcase 1.2}
\label{t5}
\end{figure}

By Lemma \ref{covering}, $F\in \F_3$. We can write $F=\{w_1uw_2\}$. If $N_G(w_1)\cap V(H)\neq \emptyset$ and $N_G(w_2)\cap V(H)= \emptyset$, we have $\F^*= (\F\backslash\{w_1uw_2\})\cup \{u'w_1u\}$ is a $K_{1,2}$-structure-cut of $G$ and $|\F^*|\leq k$ since $G-\F^*$ is disconnected and $w_2$ a component of $G-\F^*$, where $u'\in N_G(w_1)\cap V(H)$, as shown in Figure \ref{t5} (1).  The same result holds when  $N_G(w_1)\cap V(H)= \emptyset$ and $N_G(w_2)\cap V(H)\neq \emptyset$. This contradicts our assumption.

If $N_G(w_1)\cap V(H)\neq \emptyset$ and $N_G(w_2)\cap V(H)\neq \emptyset$, we have $\F^*= (\F\backslash\{w_1uw_2\})\cup \{w_1u'w_2\}$ is a $K_{1,2}$-structure-cut of $G$ and $|\F^*|\leq k$ since $G-\F^*$ is disconnected and $u$ a component of $G-\F^*$ when $N_G(w_1)\cap N_G(w_2)\cap V(H)\neq \emptyset$, where $u'\in N_G(w_1)\cap N_G(w_2)\cap V(H)$, as shown in Figure \ref{t5} (2). This contradicts our hypothesis.

If $N_G(w_1)\cap N_G(w_2)\cap V(H)=\emptyset$, then for any vertex $u'\in N_G(w_1)\cap V(H)$ and $u''\in N_G(w_2)\cap V(H)$, different neighbors of $u'$ and $u''$ can be found in $H$, labeled $v'$ and $v''$ respectively, since $H$ is 2-connected and $|V(H)|\geq 4$. We can obtain that $\F^*= (\F\backslash\{w_1uw_2\})\cup \{w_1u'v',w_2u''v''\}$ is a $K_{1,2}$-structure-cut of $G$ and $|\F^*|\leq k$ since $G-\F^*$ is disconnected and $u$ a component of $G-\F^*$ or $G-\F^*$ is trival, as shown in Figure \ref{t5} (3). This contradicts our hypothesis.

\emph{Case 2: $\F_b\neq \emptyset$.}

For any $F$ of $\F_{b}$, we write $F=\{uwv\}$, where $w\in X$ and $\{u,v\}\subseteq S_1$. If $N_G(w)\cap V(H)= \emptyset$, we have $\F^*=\F\backslash\{uwv\}$ is a $K_{1,2}$-structure-cut of $G$ and $|\F^*|<k$ since $G-\F^*$ is disconnected and $uwv$ a component of $G-\F^*$, contradicting the given condition.

If $N_G(w)\cap V(H)\neq \emptyset$, we have $\F^*=(\F\backslash\{uwv\})\cup\{uwu'\}$ is a $K_{1,2}$-structure-cut of $G$ and $|\F^*|\leq k$ since $G-\F^*$ is disconnected and $v$ a component of $G-\F^*$, where $u'\in N_G(w)\cap V(H)$, a contradiction.

\emph{Case 3: $\F_c\neq \emptyset$.}

For any $F$ of $\F_{c}$, we write $F=uwf$, where $w\in X$ and $u\in S_1$ and $f\in S_2$. If $N_G(w)\cap V(H)=N_G(f)\cap V(H)=\emptyset$, we have $\F^*=\F\backslash\{uwf\}$ is a $K_{1,2}$-structure-cut of $G$ and $|\F^*|<k$ since $G-\F^*$ is disconnected and $uwf$ a component of $G-\F^*$, contradicting the given condition.

If $N_G(w)\cap V(H)=\emptyset$ and $N_G(f)\cap V(H)\neq\emptyset$,then  we have $\F^*=(\F\backslash\{uwf\})\cup\{wff'\}$ is a $K_{1,2}$-structure-cut of $G$ and $|\F^*|\leq k$ since $G-\F^*$ is disconnected and $u$ a component of $G-\F^*$, where $f'\in N_G(f)\cap V(H)$. This is in contradiction with our assumption.

If $N_G(w)\cap V(H)\neq\emptyset$, then we have $\F^*=(\F\backslash\{uwf\})\cup\{fwf''\}$ is a $K_{1,2}$-structure-cut of $G$ and $|\F^*|\leq k$ since $G-\F^*$ is disconnected and $u$ a component of $G-\F^*$, where $f''\in N_G(w)\cap V(H)$. This is in contradiction with our assumption.

Based on the above discussion, we conclude that $\kappa(G; K_{1,2})\leq \kappa(G)$. Moreover, it is easy to see that $\kappa(C_n; K_{1,2})= \kappa(C_n)=2$ for $n\geq 7$.\endproof

\section{Conclusions}
In this paper, we have state some results on the existence of $K_{1,2}$-structure-cut in a
graph, and investigate the relationship between the connectivity $\kappa(G)$ and $K_{1,2}$-structure-connectivity
$\kappa(G;K_{1,2})$ of a graph $G$, have obtained best possible upper and lower bounds of $\kappa(G;K_{1,2})$ in terms
of $\kappa(G)$. Thus, we propose the following problem: Let $G$ be a connected graph which has $K_{1,2}$-structure-cut.
Does $\kappa(G)/m\leq\kappa(G; K_{1,m})\leq\kappa(G)$ hold for any integer $m\geq 4$? However,  However, according to the procedure of the proofs of  Lemma \ref{covering} and Theorem \ref{result}, we think it would be difficult to verify this using the same method. Some new ideas should be proposed to further pursue this
research.

\section*{Acknowledgements}
This work was supported by the National Natural Science Foundation of China (No. 12101203)
and China Postdoctoral Science Foundation (No. 2022M711075).

\end{document}